\documentclass[12pt]{amsart}

\usepackage{amsxtra,amssymb,amsthm,amsmath,amscd,url,xypic,mathrsfs}
\usepackage[latin1]{inputenc}
\usepackage{fullpage}
\usepackage{supertabular}

\usepackage[
        colorlinks,
        backref,
]{hyperref}
\usepackage{url}
\numberwithin{equation}{section}

\renewcommand{\leq}{\leqslant}
\renewcommand{\geq}{\geqslant}

\newcounter{margin}
{\end{itshape}  \bigskip}

\DeclareMathOperator{\Res}{Res}
\DeclareMathOperator{\id}{id}

\newcommand{\bfP}{\mathbb P}

\newcommand{\C}{\mathbb C}

\newcommand{\Z}{\mathbb Z}
\newcommand{\R}{\mathbb R}
\newcommand{\N}{\mathbb N}

\newcommand{\calV}{\mathcal{V}}

\newcommand{\eps}{\varepsilon}

\newcommand{\lra}{\longrightarrow}

\DeclareMathOperator{\Gal}{Gal}

\DeclareMathOperator{\im}{Im}
\DeclareMathOperator{\Tr}{Tr}
\DeclareMathOperator{\Hom}{Hom}
\DeclareMathOperator{\End}{End}

\DeclareMathOperator{\disc}{disc}

\DeclareMathOperator{\GL}{GL}

\DeclareMathSymbol{\gena}{\mathord}{letters}{"3C}
\DeclareMathSymbol{\genb}{\mathord}{letters}{"3E}

\def\rad{\operatorname{rad}}

\theoremstyle{plain}
\newtheorem{theorem}{Theorem}[section]

\newtheorem{lemma}[theorem]{Lemma}
\newtheorem{corollary}[theorem]{Corollary}

\newtheorem{proposition}[theorem]{Proposition}

\theoremstyle{remark}
\newtheorem{remark}[theorem]{Remark}

\theoremstyle{definition}

\newcommand{\nsp}{{\rm N}_{\rm spin}}

\begin{document}

\title[Computations involving Bezoutians]{On the bilinear structure
  \\associated to Bezoutians}

\author{F. Jouve}

\author{F. Rodriguez-Villegas}

\thanks{FRV was supported by NSF grant DMS-0800099. He would also like
  to thank the Centre de Recerca Matem\`atica in Barcelona for its
  hospitality}

\maketitle

\begin{abstract}
  This paper is partly a survey of known results on quadratic forms
  that are hard to find in the literature.  Our main focus is a
  twisted form of a construction due to Bezout. This {\it skew
    Bezoutian} is a symplectic (resp. quadratic) space associated to a
  pair of reciprocal (or skew reciprocal) coprime polynomials of same
  degree.  The isometry group of this space turns out to contain a
  certain associated hypergeometric group .  Using the skew Bezoutian
  we construct explicit isometries of bilinear spaces with given
  invariants (such as the characteristic polynomial or Jordan form
  and, in the quadratic case, the spinor norm).
\end{abstract}

\section{Introduction}
 \label{remark-eva}
 Given a finite dimensional commutative algebra $A$ over a field $k$
 it is useful to have a non-trivial {\it transfer map} that takes
 bilinear modules over $A$ to bilinear modules over $k$ preserving
 non-degeneracy. A typical case is that of $A=K$ a separable finite
 field extension where one can use the usual trace $\Tr_{K/k}$ for
 this purpose; for an inseparable field extension $K/k$ however, the
 trace $\Tr_{K/k}$ is identically zero.  Nevertheless, one may still
 find suitable linear maps (see~\cite[Remark $1.4$]{Mil} 
 and~\cite[discussion preceding Prop. $(2.2)$]{Ba}).

 As an example, consider the {\it unit form} $\Psi(x,y)=xy$ on a
 finite dimensional $k$-algebra (with unit) $A$. It is clearly
 non-degenerate. Transferring this form to $k$ consists of finding a
 $k$-linear map $t: A\rightarrow k$ such that the $k$-bilinear form
 $t\circ\Psi (x,y)=t(xy)$ is non-degenerate as well. The resulting
 pair $(A,t)$ is called a {\it Frobenius algebra}.

 In \S\ref{ClassicalBez} we study in some detail the situation where
 the Frobenius algebra $A$ is {\it monogenic}, i.e., $A=k[\alpha]$ for
 some $\alpha \in A$. We show that the isomorphism classes of these
 algebras over $k$ are parametrized by rational functions $w\in k[T]$
 with $w(\infty)=0$. In turn, the associated bilinear form is
 essentially given by the classical Bezoutian of the polynomials $p$
 and $q$, where $w=p/q$ in lowest terms.

 In \S\ref{skew-bezoutian} we further assume the characteristic of the
 base field $k$ to be different from $2$ and we study a skew version of the
 classical Bezoutian, which turns out to be quite interesting. For
 example, we show how it gives a natural description of the Hermitian
 form fixed by an associated hypergeometric group.

 The rest of the paper is devoted to applications
  (in characteristic different from $2$) of the skew
 Bezoutian to the problem of the existence of isometries with
 prescribed characteristic polynomial (and/or spinor norm in the
 quadratic case) or with prescribed Jordan form.

 Using the trace map and other linear maps to transfer quadratic 
 $A$-modules to quadratic
 $k$-modules, where $A$ is a separable $k$-algebra, appears
 prominently in the literature.  A general account of applications of
 trace forms to the construction of lattices via number fields can be
 found in~\cite{Bayer}.  See also~\cite[Section
 $1$]{Mil},~\cite[pp. 109--110]{KW} or~\cite[beginning of Section
 $2$]{Ba} for a general construction. 
 
 Transfer constructions are used in knot theory at least since the 1960's; a survey of some 
 of the results can be found in~\cite{KW}. The trace plays a prominent role, but other 
 transfers are also used, in particular by Trotter (see for instance~\cite[pp. 292--294]{Tr2} 
 and~\cite[pp. 181--182]{Tr1}). Moreover a general study of transfers and their applications 
 was started in the late 1960's by Scharlau~(\cite{Sch1},~\cite{Sch2}; 
 see also~\cite[Chap. 7]{Lam} for a detailed exposition of Scharlau's results), and still 
 plays a very important role in the study of quadratic and hermitian forms.

 \par
\medskip
\par
\textbf{Acknowledgements.} We would like to thank Eva Bayer for comments and suggestions 
on an earlier version of the paper as well as for pointing us to the relevant literature especially 
regarding references from knot theory.

 \section{Monogenic Frobenius algebras and the
   Bezoutian} 
\label{ClassicalBez}
Let $k$ be a field. A \emph{monogenic Frobenius algebra over $k$}, MFA
for short, is a triple $(A,\alpha,t)$, where $A$ is a finite
dimensional $k$-algebra, $A=k[\alpha]$, and $t\colon A\rightarrow k$
is a linear map such that the bilinear form
$$
\langle x,y\rangle:=t(xy)\,,\qquad x,y\in A\,,
$$
is non-degenerate.  With some notation abuse we will say in this case
that the linear map $t$ is non-degenerate.  Two such algebras
$(A,\alpha,t),(A',\alpha',t')$ are isomorphic if there is an
isomorphism of algebras $\phi: A\rightarrow A'$ such that
$\alpha'=\phi(\alpha)$ and $t'=t\circ\phi^{-1}$.

Let $A$ be a $k$-algebra and $t: A \rightarrow k$ be a non-degenerate
linear map.  For $a\in A$ the map $t_a(x):=t(ax)$ is also
non-degenerate if and only if $a$ is not a zero-divisor. In
particular, there is natural action of the group of units $A^\times$
of $A$ on the set of non-degenerate linear maps $t: A \rightarrow k$
defined by setting $(a\cdot t)(x):=t_a(x)=t(ax)$, where $a\in
A^\times$.

The following theorem gives a parametrization of MFA's.
\begin{theorem}\label{param-transfer-alg}
  For any $d\geq 1$ the map $(A,\alpha,t)\mapsto w(T)$, where
$$
w(T):=\sum_{\ell\geq 0}t(\alpha^\ell)T^{-\ell-1}\in k[[T^{-1}]],
$$ 
induces a bijection between isomorphism classes of $d$-dimensional
MFA's $(A,\alpha,t)$ over $k$ and rational functions $w\in
k(T)\cap k[[T^{-1}]]$ of degree $d$ with $w(\infty)=0$.
\end{theorem}

\begin{proof}
  Let $(A,\alpha,t)$ be an $d$-dimensional MFA over $k$.
  Consider the power series
$$
w(T):=\sum_{\ell\geq 0}t(\alpha^\ell)T^{-\ell-1}\in k[[T^{-1}]]\,.
$$ 
Let $q=T^d+\sum_{i=0}^{d-1}q_iT^i\in k[T]$ be the minimal polynomial
of $\alpha$ over $k$. Set $q_d=1$ and $q_i=0$ for $i>d$. We have:
\begin{equation}
\label{qw-prod}
q(T)w(T)=\sum_{\ell\in \Z} \sum_{i\geq
  \ell+1}q_it(\alpha^{i-\ell-1})\,T^\ell\,. 
\end{equation}
On the right hand side the coefficient of $T^\ell$ vanishes as soon as
$\ell\geq d$ since $q_i=0$ for $i>d$. Moreover, since $q$ vanishes at
$\alpha$,  we have for any integer $\ell$ 
$$
0=\alpha^{-\ell-1}\sum_{i\geq 0}q_i\alpha^i=\sum_{i\geq 0}q_i\alpha^{i-\ell-1}\,.
$$
Thus the coefficient of $T^\ell$ on the right hand side
of~\eqref{qw-prod} also vanishes for $\ell+1\leq 0$. We deduce that
$q(T)w(T)$ is a polynomial $p\in k[T]$ of degree $\leq d-1$ with
$$
p(T):=\sum_{\ell=0}^{d-1}
\sum_{i=\ell+1}^dq_it(\alpha^{i-\ell-1})\,T^\ell\,.  
$$
The rational function $w=p/q\in k(T)$ satisfies $w(\infty)=0$ since
$\deg(p)<\deg(q)$. This calculation is valid for any linear map $t$
without assuming it is non-degenerate.

 The Gram matrix of
$\langle\cdot,\cdot\rangle$ in the $k$-basis
$(1,\alpha,\ldots,\alpha^{d-1})$ of $A$ is the Hankel matrix
\begin{equation}
\label{hankel}
H(p/q):=(t(\alpha^{i+j}))_{0\leq i,j\leq d-1}.
\end{equation}
By assumption the bilinear form $\langle\cdot,\cdot\rangle$ is
non-degenerate hence the determinant of $H(p/q)$ is non-zero. By
Kronecker's theorem (in~\cite{Kr}, see e.g.~\cite[Th. 8.20 \& Prop. 8.22]{Fu}) $p$ and $q$ are coprime and $\deg
w=d$.  Obviously $w$ depends only on the isomorphism
class of $(A,\alpha, t)$. Indeed if $(A,\alpha,t)$ and $(A',\alpha',t')$ are two 
isomorphic MFA's and if $\phi$ is a fixed isomorphism then for all $\ell\geq 0$,
$$
t'(\alpha'^\ell)=t\circ \phi^{-1}\left(\phi(\alpha)^\ell\right)=t(\alpha^\ell)\,.
$$

 Conversely, suppose we are given a rational function
$w=p(T)/q(T)\in k(T)$ satisfying $w(\infty)=0$ and $\deg(w)=d\geq
1$.   Without loss of generality we
may assume $q$ monic, $\deg q=d$, $(p,q)=1$, $p$ non-zero and $\deg
p<d$. Set
$$
A:=k[T]/(q)\,.
$$
Let $\alpha$ be the image of $T$ in $A$; then
$\{1,\alpha,\ldots,\alpha^{d-1}\}$ is a $k$-basis of $A$.  Let
$\sum_{\ell\geq 0}t_\ell T^{-\ell-1}\in k[[T^{-1}]]$ be the power
series expansion of $w$ at $\infty$.  Define the $k$-linear map
\begin{equation}
\label{t-defn}
t\colon A\rightarrow k\,,\qquad \alpha^\ell\mapsto t_\ell, \quad
\ell=0,1,\ldots, d-1\,.
\end{equation}
By construction the power series in $T^{-1}$:
$$
\tilde{w}(T):=\sum_{\ell\geq 0}t(\alpha^\ell)T^{-\ell-1}\,,
$$
has the same first $d$ coefficients as $w$. On the other hand, as
observed above, $\tilde p(T) := q(T)\tilde{w}(T)$ is a polynomial of
degree $\leq d-1$.  Hence the first $d$ coefficients of $\tilde p(T)$ and $q(T)w(T)= p(T)$ agree
and it follows that $w=\tilde{w}$.

Again, since $(p,q)=1$, by Kronecker's theorem the matrix $H(p/q)$ has
non-zero determinant and hence the bilinear form on $A$ defined by
$(x,y)\mapsto t(xy)$ is non-degenerate. This completes the proof of
the theorem.
\end{proof}

\begin{corollary}
\label{main-corollary}
1) Given a monic polynomial $q\in k[T]$ of degree $\deg(q)>0$ there
exists a MFA of the form $(k[T]/(q),T \bmod q,t)$.

2) For any MFA $(A,\alpha,t)$ the unit group $A^\times$ acts
transitively on the non-degenerate linear maps on $A$.
\end{corollary}
\begin{proof}
  To prove 1), it is enough to take $w=1/q$. Let $(k[T]/(q), T \bmod
  q,t)$ be the corresponding MFA. Say $t'$ is another non-degenerate
  linear map on it and let $p/q$ be the associated rational
  function. To show 2), it is enough to prove that
$$
\sum_{l\geq 0} t(p(\alpha)\alpha^l)\,T^{-l-1}= \frac{p(T)}{q(T)},
$$
where $\alpha:= T\bmod q$, since then $t'=p(\alpha)\cdot t$.  To see
this note that if $p(T)=\sum_{j=0}^{d-1}p_j\,T^j$ then the left hand
side equals
$$
\sum_{l\geq 0}\sum_{j=0}^{d-1} p_jt_{l+j}\,T^{-l-1}= 
\sum_{l\geq 0}\sum_{j=0}^{d-1} p_jt_l\,T^{j-l-1}= p(T)\cdot \frac 1{q(T)}.
$$
finishing the proof.
\end{proof}

\begin{remark}
  Over the complex numbers the space of rational functions $w$ of
  degree $d$ with $w(\infty)=0$ is naturally isomorphic to a circle
  bundle over the moduli space of $SU(2)$ monopoles of charge
  $d$~\cite{Do}. 
\end{remark}

\subsection{Examples}
1) If $q$ is irreducible and separable (i.e., $q$ is irreducible and
$k$ has characteristic $0$ or $q$ is irreducible, $k$ has
characteristic $l>0$, and $q$ is not a polynomial in $T^l$) then
$A:=k[T]/(q)$ is a field, $K$ say, and it is well-known that
$t=\Tr_{K/k}$ is non-degenerate (in fact the algebra $A$ is separable
over $k$ if and only if $\Tr_{A/k}$ is non-degenerate).  It is not
hard to see that the underlying MFA corresponds to the rational
function $p/q$, where $p\equiv dq/dT \bmod q$. Indeed let $L$ be the
Galois closure of $K/k$. Since the extension $K/k$ is separable there
are $d:=[K:k]$ distinct $k$-embeddings $\sigma_1,\ldots ,\sigma_d:\,
K\hookrightarrow L$. The Galois action of $G:=\Gal(L/K)$ on $L$
extends to an action on $L(T)$ via $\sigma(\sum
\lambda_iT^i):=\sum\sigma(\lambda_i)T^i$. Therefore
$$
w=\sum_{i\geq 0}\Tr_{K/k}(\alpha^i)T^{-i-1}=\sum_{i\geq
  0}\left(\sum_{j=1}^d\sigma_j(\alpha^i)\right)T^{-i-1} ={1\over
  T}\sum_{j=1}^d\sigma_j\left(\sum_{i\geq 0}(\alpha/T)^i\right)\,.
$$
The inner sum equals $T/(T-\alpha)$ hence:
$$
w=\sum_{j=1}^d\sigma_j\left({1\over
    T-\alpha}\right)=\sum_{j=1}^d{1\over
  T-\sigma_j(\alpha)}={dq/dT\over q}\,. 
$$

2) Another extreme case is $A=k[T]/(T^d)$. If say $p=1$ then the
rational function $w$ given by Theorem~\ref{param-transfer-alg} is
simply $1/T^d$. The $k$-linear map $t: k[T]/(T^d)\rightarrow k$
corresponding to $w$ is defined by $t(\alpha^i)=0$ for $0\leq i\leq
d-2$ and $t(\alpha^{d-1})=1$.  Thus $t$ can be identified with the
projection $A\rightarrow A$ with image $k\alpha^{d-1}$. The
non-degeneracy of $t$ can be shown by elementary arguments. Namely, if
$z=\sum z_i\alpha^i\in A$ is orthogonal to any $y\in A$ with
respect to the inner product $(x,y)\mapsto t(xy)$ then in particular,
for any fixed index $i$ we have
$t(z\cdot\alpha^{d-1-i})=z_i=0$. Thus $z=0$.

\subsection{Reproducing kernel}
Given a Frobenius algebra $(A,t)$ of dimension $d$ consider the
\emph{Casimir element} (or \emph{reproducing kernel}) defined by
$$
C:=\sum_{i=1}^d e_i \otimes e_i^\# \in A\otimes A,
$$
where $e_1,\ldots,e_d$ is any basis of $A$ over $k$ and $e_i^\#$ is
its dual basis (with respect to the bilinear form $\langle
\cdot,\cdot\rangle$ determined by $t$), i.e.,
$$
\langle e_i,e_j^\#\rangle=t(e_ie_j^\#)=\delta_{i,j}, \qquad
i,j=1,\ldots, d.
$$
This element is well defined; it does not depend on the choice of
basis used in its definition.

We have
$$
C=\sum_{i,j=1}^d \langle e_i^\#,e_j^\#\rangle\, e_i\otimes e_j.
$$
For a MFA $(A,\alpha,t)$ with $A\simeq k[T]/(q)$ for some $q\in k[T]$
monic of degree $d$ (namely, the minimal polynomial of $\alpha$) we
can represent elements of $A\otimes A$ as polynomials in $k[x,y]$ of
degree at most $d-1$ in each variable.

 Taking $e_i:=\alpha^{i-1}$ for $i=1,\ldots, d$ as our basis of $A$ we
 obtain 
$$
C=\sum_{i,j=1}^d b^\#_{i,j} \,x^{i-1}y^{j-1},
$$
where $b^\#_{i,j}:=\langle e_i^\#,e_j^\#\rangle$. 

The matrix $B^\#:=(b^\#_{i,j})$ is the inverse of the Hankel
matrix $H(p/q)$ in~\eqref{hankel} since the matrices 
$$
\left(\langle e_i,e_j\rangle\right), \qquad \left(\langle
  e_i^\#,e_j^\#\rangle\right) 
$$
are inverses of each other.  Combining this observation with the proof
of Corollary~\ref{main-corollary} we obtain the following.
\begin{proposition}
\label{prop1}
With the above notation and assumptions
$$
H(1/q)M_p={}^tM_pH(1/q)=H(p/q)=(B^\#)^{-1},
$$
where $M_p$ is the matrix of multiplication by $p$ in $k[T]/(q)$ in
the basis $1,T,\ldots,T^{d-1}$.
\end{proposition}

\subsection{Classical Bezoutian}

Given two polynomials $p,q\in k[T]$, the {\it classical Bezoutian}
is the  symmetric matrix $B(p,q):=(b_{i,j})$, where 
\begin{equation}
\label{bez-defn}
\frac{p(x)q(y)-p(y)q(x)}{x-y}=\sum_{i,j=1}^d b_{i,j}\,x^{i-1}y^{j-1},
\end{equation}
with $d:=\max\{\deg(p),\deg(q)\}$.
We have the following matrix expression for $B(p,q)$.
\begin{lemma}
\label{prop2}
Let $Q$ be the matrix with entries $(q_{i+j-1})$, where $1\leq i,j\leq d$ and
$q=\sum_{i\geq 0}q_i\,T^i\in k[T]$ is a polynomial of degree $d$. Assume $\deg p\leq d$, 
then
$$
B(p,q)=-M_pQ,
$$
where $M_p$ is the matrix of multiplication by $p$ in $k[T]/(q)$.
\end{lemma}
\begin{proof}
  Expand the left hand side of~\eqref{bez-defn} in Laurent series in
  $k[x,y][[y^{-1}]]$. Modulo $q(x)$ the coefficient of $y^{j-1}$ is that of
$$
-p(x)q(y)y^{-1}\sum_{i\geq 0} \left(\frac x y\right)^i
$$
and this is easily seen to equal $-p(x)\sum_{i\geq 0}q_{i+j}x^i.$
\end{proof}

We leave the easy proof of the following lemma to the reader
\begin{lemma}
\label{Q-lemma}
We have
$$
Q=H(1/q)^{-1}.
$$
\end{lemma}

Putting together Propositions~\ref{prop1} and Lemmas~\ref{prop2}
and~\ref{Q-lemma} we finally obtain a connection between the Casimir
element $C$ and the classical Bezoutian.
\begin{theorem}\label{thBsharp}
With the above notation and assumptions
$$
B^\#=B(q,r),
$$
where $r\in k[T]$ is a polynomial of degree less than $d$ such that
$rp\equiv 1\bmod q$. Or, equivalently,
$$
C=\frac{q(x)r(y)-q(y)r(x)}{x-y}.
$$
\end{theorem}

\begin{remark} We have discussed two symmetric matrices associated to a
  pair of coprime polynomials $p,q\in k[T]$: the Hankel matrix
  $H(p/q)$ and the Bezoutian $B(p,q)$. Combining
  Propositions~\ref{prop1},~\ref{prop2} and Lemma~\ref{Q-lemma} we
  find that they are actually congruent up to a minus sign
$$
{}^tQH(p/q)Q= -B(p,q),
$$
(note that $Q$ is symmetric).
\end{remark}

The Bezoutian plays an important role in mathematical control theory,
see for example~\cite{Fu}.

\section{The skew Bezoutian}
\label{skew-bezoutian}

We now turn to a construction that is a skew version of  the
classical Bezoutian.

\subsection{Preliminaries}
\label{prelim}

Let $k$ be a field and $d$ a positive integer; consider the algebra
$A:=k[T]/(T^d)$.  Given a power series $a=a_0+a_1T+a_2T^2+\cdots \in
k[[T]]$ let $M(a) \in k^{d\times d}$ be the matrix of the
$k$-linear map $A\rightarrow A$ defined by multiplication by $a$ in
the basis $1,T,\ldots,T^{d-1}$ of $A$. Concretely,
 $$
M(a) =
\left(
\begin{array}{cccc}
a_0 & 0 &  & 0 \\
a_1 & a_0 &  & \\
\vdots &  & \ddots & 0 \\
a_{d-1} & \cdots & \cdots & a_0
\end{array}
\right).
$$
The map $a \mapsto M(a)$ is clearly a homomorphism of $k$-algebras
$k[[T]] \rightarrow M^{d\times d}(k)$.

We will assume from now on that $k$ has characteristic different from
$2$. For a polynomial $q\in k[T]$ we let $q^*$ be the polynomial $q$
with its coefficients reversed, i.e.,
$$
q^*(T):=T^{\deg(q)}q(1/T), \qquad q\in k[T],
$$
where $\deg(q)$ is the degree of $q$.

For further reference let us note a few simple observations about
the operation $*$. In general, $*$ is not additive but we have
$$
(p+q)^*=p^* +q^*, \quad \quad \text{if } \deg(p)=\deg(q)
$$
and $(pq)^*=p^* q^*$ always. We extend $*$ to $k(T)$ by
multiplicativity. Then for $w=p/q$ we have
$$
w^*(T):=\frac{p^*}{q^*}=T^{\deg(p)-\deg(q)}w(T^{-1}).
$$
We will say that $w\in k(T)$ is {\it reciprocal} if $w^*=w$, {\it
  skew-reciprocal} if $w^*=-w$ and in general {\it $\eps$-reciprocal}
if $w^*=\eps w$ with $\eps=\pm 1$.

To shorten the notation we let $v_\eps$, for $\eps =\pm 1$,
denote the valuation on $k[T]$ at $T-\eps$. If $p\in k[T]$ is
skew-reciprocal then $p(1)=0$. It follows that in fact $v_+(p)$ must
be odd since otherwise $p(T)/(T-1)^{v_+(p)}$ would be a skew-reciprocal
polynomial not vanishing at $T=1$. Similarly, if $p\in k[T]$ is
$\eps$-reciprocal with $\eps=-(-1)^{\deg(p)}$ then $p(-1)=0$ and
again, $v_-(p)$ must be odd. In particular, if $p(-1)\neq 0$ then
$\deg(p)$ must be even.

\subsection{Definition}

Let $w\in k(T)$ be a rational function with coefficients in $k$.
Assume that $w$ is regular at $0$ and $\infty$. Then we have the two
power series expansions
$$
w(T)=w_0+w_1T +\cdots, \qquad
w(T)=w_0^*+w_1^*T^{-1} +\cdots
$$
Let $d$ be the degree of $w=p/q$ where the fraction is written in lowest terms.
 With the above notation define the {\it
  skew Bezoutian} of $w$ as
$$
B^*(w)={}^tM(T^{d-\deg p}w^*)-M(w)
=
\left(
\begin{array}{cccc}
w_0^*-w_0 & w_1^* & \ldots & w_{d-1}^* \\
-w_1 & w_0^*-w_0 & w_1^* & w_{d-2}^*\\
\vdots & \ddots & \ddots & w_1^* \\
-w_{d-1} & \cdots & -w_1 & w_0^*-w_0
\end{array}
\right).
$$
(We learned of this construction in~\cite{HA}). Note that $B^*(w)$ is
a {\it Toeplitz} matrix (constant entries along
diagonals). Also, in the case where $d=\deg p=\deg q$ (for instance if $w(\infty)=1$, as will be assumed later) one has the simpler definition 
$B^*(w)={}^tM(w^*)-M(w)$.

There is a more conceptual way to give $B^*$, closer to the approach
of the previous section on Frobenius algebras, as follows.  Let
$R:=k[T,T^{-1}]$ and let $t:R\rightarrow k$ be the linear map
corresponding to taking constant terms; i.e., $t(1):=1$ and
$t(T^n):=0$ for all non-zero integers $n$.

We may represent elements in the dual space $\omega\in \Hom(R,k)$ as
formal infinite series of the form
$$
\omega:=\sum_{n\in \Z}\omega(T^n)\,T^{-n}.
$$
Following the usual rules of multiplication of series gives
$\Hom(R,k)$ the structure of an $R$-module. Then
$$
\omega(u)=t(u\cdot \omega), \qquad u \in R.
$$

With this notation define $\omega\in \Hom(R,k)$ by
$$
\omega:=\sum_{n\geq 0}w_n^*T^n - \sum_{n\geq 0}w_nT^{-n}.
$$
Explicitly,
\begin{equation}
\label{omega-defn}
\omega(T^n)=
\begin{cases}
-w_n & n>0\\
w_0^*-w_0 & n=0\\
w_{-n}^* & n <0
\end{cases}.
\end{equation}
Then the skew Bezoutian $B^*$ is the matrix with entries $\omega(T^{i-j})$ for
$i,j=0,1,\ldots, d-1$.

The first thing to point out is the value of the
determinant of $B^*(w)$. Write $w=p/q$ for  $p,q\in k[T]$ relatively
prime. By assumption $d=\deg(w)=\max\{\deg(p),\deg(q)\}$. But since we
also assume $w$ is regular at infinity we must have $\deg(p)\leq
\deg(q)=d$. 
\begin{proposition}
\label{bez-res}
Let $q_0$ and $q_d$ be the constant and leading coefficients of $q$
respectively. Then we have
$$
q_0^dq_d^{\deg(p)}\det B^*(w)=(-1)^{d-\deg(p)}\Res(p,q)
$$
\end{proposition}
\begin{proof}
  Assume first that $\deg(p)=d$.  The following block-matrix identity
  is easy to check using the fact that $M$ is a homomorphism. 
$$
\left(
\begin{array}{cc}
{}^tM(w^*) &  M(w)\\
I_d & I_d  
\end{array}
\right)
\cdot
\left(
\begin{array}{cc}
{}^tM(q^*) & 0 \\
0 &  M(q)
\end{array}
\right)
=
\left(
\begin{array}{cc}
{}^tM(p^*) & M(p)  \\
{}^tM(q^*) &  M(q)
\end{array}
\right)\,,
$$
where $I_d$ is the $d\times d$ identity matrix.  The right hand side
is then precisely the Sylvester matrix of $p$ and $q$ whose
determinant is $\Res(p,q)$.  On the other hand, the determinant of the
left hand side equals $(q_0q_d)^d\det B^*(w)$.

We may consider the case $\deg(p)<d$ as a specialization of the
generic case of $\deg(p)=d$.  Then the determinant on the right hand
side is easily seen to equal $(-q_d)^{d-\deg(p)}\Res(p,q)$ completing the
proof.
\end{proof}
\begin{remark}
  The proof we gave follows that of a generalization
  of Proposition~\ref{bez-res} to subresultants in~\cite[Prop. 11]{BSP}. The
  fact that the classical Bezoutian has determinant related to the
  resultant goes back to Bezout. From a computational point of view
  the Bezoutian has the advantage that it is a matrix of size
  $\max\{\deg(p),\deg(q)\}$ versus the Sylvester matrix that has size
  $\deg(p)+\deg(q)$.
\end{remark}

Consider the case that $\deg(p)=d$. It follows from the proof of the
proposition that
\begin{equation}
\label{skew-bez-pol}
B^*(p,q):={}^tM(q^*)B^*(w)M(q)={}^tM(p^*)M(q)-{}^tM(q^*)M(p)
\end{equation}
has determinant $\Res(p,q)$. This matrix can be described in a very
similar way to that of the classical
Bezoutian~\eqref{bez-defn}. Indeed, it is not hard to see that its
entries are the coefficients of the two variable polynomial
\begin{equation}
\label{skew-bez-pol-1}
\frac{p(x)q^*(y)-q(x)p^*(y)}{xy-1}.
\end{equation}

As in the proof of the proposition we may think of the case
$\deg(p)<d$ as a specialization of the generic case $\deg(p)=d$ and
define $B^*(p,q)$ accordingly (namely, $p^*$ should be replaced by
$T^dp(T^{-1})$). In general we have
$$
\det(B^*(p,q))=(-q_d)^{d-\deg(p)}\Res(p,q).
$$

\subsection{Bilinear form}

We now consider the bilinear form determined by the skew
Bezoutian. Let $V:=R/(q^*)$ with basis $1,T,\ldots, T^{d-1}$. We claim
that the linear form $\omega$~\eqref{omega-defn} vanishes on the ideal
$(q^*)\subseteq R$ and therefore induces a corresponding linear form
on $V$. To see this we compute
$$
q^*\cdot \omega= q^*\sum_{n\geq 0} w_n^* T^n- q^*\sum_{n\geq 0}
w_nT^{-n}.
$$
The first term equals
$q^*(T)w(T^{-1})=T^{d-\deg(p)}q^*(T)w^*(T)=T^{d-\deg(p)}p^*(T)$ and
similarly, the second term also equals $q^*(T)w(T^{-1})$ cancelling
out. It follows that $\omega(T^nq^*(T))=t(T^nq^*(T)\cdot \omega)=0$
for all $n$.

{}From now on we assume 
\begin{equation}
\label{w-fctn-eqn}
  w(T^{-1})=-\eps w(T), \qquad w(\infty)=1,
\end{equation}
for some $\eps=\pm1$.
Then
$$
{}^tB^*(w)=\eps B^*(w).
$$
and, in the notation of the previous section,
$$
w_n=-\eps w_n^*, \qquad n=0,1,\ldots.
$$
Therefore, 
\begin{equation}
\label{sbez-defn}
B^*(w)=
\left(
\begin{array}{cccc}
\eps+1 & w_1^* & \cdots & w_{d-1}^* \\
\eps w_1^* & \eps+1 & \cdots & w_{d-2}^*\\
\vdots & \vdots & \vdots & \vdots \\
\eps w_{d-1}^* & \eps w_{d-2}^* & \cdots & \eps+1
\end{array}
\right)\,.
\end{equation}

Recall that $w=p/q$ is in lowest terms. The
assumptions~\eqref{w-fctn-eqn} imply that both $p$ and $q$ must be
reciprocal or skew-reciprocal polynomials of degree $d$. Let $\eps_p$
and $\eps_q$ be the corresponding signs: $p^*=\eps_pp$ and
$q^*=\eps_qq$ with $\eps=\eps_p\eps_q$. Note that we cannot have
$\eps_p=\eps_q=-1$ since $p$ and $q$ are relatively prime by
assumption.

We may give the {\it skew Bezoutian} bilinear form in a way analogous
to the classical one. The new ingredient is the involution
$\iota:T\mapsto T^{-1}$ of $R$ that descends to $V$ since it fixes the
ideal $(q^*)$.  Indeed we have that $B^*(w)$ is the Gram matrix in the
basis $1,T,\ldots, T^{d-1}$ of the bilinear form $\Psi$ on $V$ defined
by
$$
\Psi(u,v):=t(uv^\iota\cdot w)=\omega(uv^\iota),
$$
where, with some notation abuse, $u,v\in V$. Concretely,
\begin{equation}
\label{toeplitz}
\Psi(T^i,T^j)=\omega(T^{i-j}), \qquad i,j \in \Z.
\end{equation}

This bilinear form satisfies
\begin{equation}
  \label{bil-symmetry}
\Psi(v,u)=\eps \Psi(u,v).
\end{equation}
We will say that $(V,\Psi)$ is an {\it $\eps$-symmetric} bilinear
space over $k$. By Proposition~\ref{bez-res} this space is
non-degenerate since $p$ and $q$ are relatively prime.

\subsection{Properties}
   
In addition to the bilinear form $\Psi$ the skew Bezoutian carries
some extra structures not shared by the classical Bezoutian. It has a
distinguished vector $v_0$, the class of the polynomial $1\in R$, with
$\Psi(v_0,v_0)=1+\eps$ and an isometry $\gamma$, given by
multiplication by $T$ (the fact that it is an isometry is clearly seen
in \eqref{toeplitz}, for example). Note that by construction $\gamma$
has characteristic polynomial $\pm q$ (the monic generator of the ideal
$(q^*)=(q)$). Moreover, the translates
$v_0,\gamma(v_0),\gamma^2(v_0),\cdots$ generate the whole space
$V$. In fact, these properties characterize the skew Bezoutian as we
now show.

Given $v_0\in V$ with $\Psi(v_0,v_0)=1+\eps$ we define its associated
{\it $\eps$-reflection} to be the isometry given by
\begin{equation}
\label{reflection-defn}
\sigma(v):=v-\Psi(v_0,v)\,v_0.
\end{equation}
(In the skew-symmetric case $\sigma$ is usually called a {\it
  transvection}.)  Note that $\sigma$ is
of order two if $\eps=+1$ but of infinite order if $\eps=-1$. In fact,
$$
\sigma^{-1}(v):=v-\eps\Psi(v_0,v)\,v_0.
$$
We have
$$
\sigma(v_0)=-\eps v_0, \qquad \sigma(v)=v, \quad\text{if }
\Psi(v_0,v)=0.
$$
Hence $\sigma$ fixes a codimension $1$ subspace of $V$ and
$\det(\sigma)=-\eps$.

Recall that an isometry of a non-degenerate bilinear space has a
characteristic polynomial which is reciprocal or skew-reciprocal.

\begin{theorem} \label{characterization} Let $(V,\Psi)$ be a
  non-degenerate, finite dimensional, $\eps$-symmetric bilinear space
  over $k$. Suppose there exists an isometry $\gamma$ of this space
  and a vector $v_0\in V$ such that 

  (i) $\Psi(v_0,v_0)=1+\eps$

(ii) $V$ is generated by $v_0,\gamma v_0, \gamma^2 v_0, \cdots$.

\medskip Then $(V,\Psi)$ is the skew Bezoutian $B^*(w)$ with $w=
p/q$, where $q$ is the characteristic polynomial of $\gamma$ and $p$
is the characteristic polynomial of $\gamma\sigma$ with $\sigma$ the
$\eps$-reflection associated to $v_0$.
\end{theorem}
\begin{proof}
  Let $d$ be the dimension of $V$. Note that $\calV:=\{v_0,\gamma
  v_0,\ldots, \gamma^{d-1}v_0\}$ is a basis of $V$. Indeed, by
  the Cayley--Hamilton theorem $\gamma^nv_0$ is in the span of
  $\calV$   and hence so is every $v$ in $V$ by hypothesis (ii). Again by
  hypothesis (ii), $\calV$ is linearly independent.

Define for every $n\in \Z$
\begin{equation}
\label{c_n-defn}
c_n:=\Psi(\gamma^nv_0,v_0).
\end{equation}
Note that $c_{-n}=\eps c_n$.
We claim that 
$$
1+\sum_{n\geq 1}c_n\,T^n
$$
is the power series expansion of a rational function of denominator
$q$. Write $q=\sum_{k\geq 0}q_k\,T^k$. By assumption $q_{d-k}=q_k$
for $k=0,\ldots,d$ and $q_0=q_d=1$. Then
$$
q(T)(1+\sum_{n\geq 1}c_n\,T^n)=\sum_{n\geq 0}r_n\,T^n=1+\sum_{n\geq 1}r_n\,T^n,
$$
where $r_n=q_n+\sum_{k=1}^{n}c_k\,q_{n-k}$ for $n\geqslant 1$. Since $q_n=0$ for $n>d$ we
have 
$$
r_n=\sum_{k=0}^dc_{n-k}\,q_k=\sum_{k=0}^dc_{n-d+k}\,q_{d-k}=\sum_{k=0}^dc_{n-d+k}\,q_k, \qquad n>d.
$$
Hence
$$
r_n=\Psi(\gamma^{n-d}q(\gamma)v_0,v_0)=0, \qquad n>d.
$$
We now show that $r_{d-n}=-\eps r_n$ for $n=0,\ldots,d$. 
 Since $q(\gamma)=0$ we have for $n < d$
$$
r_{d-n}=q_n+\sum_{k=1}^{d-n}q_{n+k}\,c_k=q_n-
\Psi\left(s_n(\gamma)\gamma^{-n}
  v_0,v_0\right),
$$
where $s_n:=\sum_{k=0}^nq_kT^k$. Hence
$$
r_{d-n}=q_n-\sum_{k=0}^nq_k\,c_{k-n}=q_n-(1+\eps)q_n
-\eps\sum_{k=0}^{n-1}q_k\,c_{n-k} = -\eps r_n.
$$
We have shown then that $p(T):=-\eps\sum_{n=0}^dr_n\,T^n$ is
$(-\eps)$-reciprocal; since $r_0=1$ it is also monic. In other words,
we have that $(V,\Psi)$ is isometric to the skew Bezoutian $B^*(p,q)$.

It remains to show that $p$ is the characteristic polynomial of
$\delta:=\gamma\sigma$. For every $n\in \Z$ let $\sigma_n$ be the
$\eps$-reflection associated to $v_n:=\gamma^n v_0$. Note that
$\Psi(v_n,v_n)=1+\eps$.  We have
$$
\sigma_n=\gamma^n\sigma\gamma^{-n}
$$
and hence by induction 
$$
\delta^n=\sigma_1\cdots\sigma_n\gamma^n.
$$
Let $u_0:=v_0$ and $u_n:=\sigma_n^{-1}\cdots\sigma_1^{-1}v_0$ for
$n>0$. Let also $e_n:=\eps \Psi(\delta^nv_0,v_0)$ for $n\in \Z$. Then
$$
 e_{n+1}=\Psi(v_0,\sigma_1\cdots\sigma_{n+1}v_{n+1})=\Psi(u_{n+1},v_{n+1}). 
$$
Since
$$
u_{n+1}=\sigma_{n+1}^{-1}u_n=u_n-\eps\Psi(u_n,v_{n+1})v_{n+1}, \qquad
n\geq 0
$$ 
we get
$$
e_{n+1}=\Psi(u_n,v_{n+1})-\eps\Psi(u_n,v_{n+1})\Psi(v_{n+1},v_{n+1})=
-\eps\Psi(u_n,v_{n+1}).
$$
Therefore $u_{n+1}=u_n+e_{n+1}v_{n+1}$ and by induction
$$
u_n=v_0+\sum_{k=1}^ne_k\,v_k.
$$
Finally,
$$
-e_{n+1}=\Psi(v_{n+1},u_n)=c_{n+1}+\sum_{k=1}^ne_k\,c_{n+1-k}
$$
and
\begin{equation}
\label{gf-symmetry}
(1+\sum_{n\geq 1}c_n\,T^n)(1+\sum_{n\geq 1}e_n\,T^n)=1.
\end{equation}
Combined with our previous calculation we see that
$$
p(T)(1+\sum_{n\geq 1}e_n\,T^n)=-\eps q(T).
$$
So if $p(T)=\sum_{n=0}^d p_n\,T^n$ then
$$
0=\sum_{k=0}^de_{n-d+k}\,p_{d-k}=-\eps \sum_{k=0}^de_{n-d+k}\,p_k, \qquad n>d
$$
and
\begin{equation}
\label{p-relations}
\Psi(p(\delta)v_0,\delta^{d-n}v_0)=0, \qquad n>d.
\end{equation}
It is not hard to see that $\delta^nv_0=-\eps
v_n+\sum_{j=1}^{n-1}\alpha_{n,j}v_j$ for $n=1,\ldots,d-1$, for some
$\alpha_{n,j}\in k$ (note that for $n=1$ the equality is $\delta v_0=-\eps v_1$). It follows that the $\delta^nv_0$ with $n\in \Z$
span $V$ and by \eqref{p-relations}  $p(\delta)=0$.
\end{proof}

\begin{remark}
  The equivalence established by Theorem~\ref{characterization} is a
  skew analogue of that in Theorem~\ref{param-transfer-alg}, where the
  new ingredient is the involution $T\mapsto T^{-1}$ of the algebra
  $k[T,T^{-1}]/(q)$.  A similar result appears
  in~\cite[Prop. 3.2]{Do}.
\end{remark}

\subsection{Hypergeometric groups}

We choose now $k=\C$.  The subgroup $\Gamma \subseteq \GL(V)$
generated by $\gamma,\delta,\sigma$ (see
Theorem~\ref{characterization}) is a hypergeometric group in the sense
of~\cite[Def. 3.1]{BH}, with parameters the multisets of roots of $p$
and $q$.  In other words, we have a triple of elements in $\GL(V)$
which multiply to the identity, two of which have a prescribed
characteristic polynomial and the third fixes a codimension one
subspace of $V$. By a theorem of Levelt such triples are unique up to
conjugation by $\GL(V)$ (see~\cite[Th. 3.5]{BH}).

Since our polynomials are coprime it is proved in~\cite[Prop. 3.3 \&
Th. 4.3]{BH} that $\Gamma$ acts irreducibly.  Furthermore, if we
assume $p$ and $q$ have real coefficients, since they are $(\pm
1)$-reciprocal, $\Gamma$ fixes a non-degenerate bilinear form $\Psi$
on $V$ which is unique up to scaling.  Our discussion here shows that
this form is none other than the skew Bezoutian $B^*(p,q)$ of $p$ and
$q$.  This was mentioned in~\cite{RV}. For general $p$ and $q$ with
complex coefficients it is not hard to extend the construction of the
skew-Bezoutian and this now yields a Hermitian form fixed by
$\Gamma$.

Over $\R$ the signature $\sigma$ of $\Psi$ can be computed by a skew
version of the classical theorem of Hermite~\cite[p. 409]{He} for the usual
Bezoutian. (Here $\varsigma:=r-s$ if $\Psi$ is isometric to $x_1^2+\cdots
x_r^2-y_1^2-\cdots -y_s^2$ over $\R$.) For the classical Bezoutian the
signature depends on the interlacing pattern of the roots of $p$ and
$q$ in $\R$.  For the skew Bezoutian it depends on the interlacing
pattern of the roots on the unit circle $S^1$. In both cases this can
be phrased in terms of the Cauchy index for the rational function
$p/q$ (on $\R$ for the classical case, on $S^1$ for the skew case; for
the latter see~\cite[Th. 2.1]{HA}).

A conceptual formulation of Hermite's result is as follows. The
rational function $w=p/q\in \R[T]$ gives a continuous map
$w:\bfP^1(\R)\rightarrow \bfP^1(\R)$. In turn this yields a
homomorphism $H_1(\bfP^1(\R),\Z)\rightarrow H_1(\bfP^1(\R),\Z)$. After
fixing an isomorphism $H_1(\bfP^1(\R),\Z)\simeq \Z$ this map is
multiplication by some integer which is none other than the signature
$\varsigma$.  The same applies for the skew Bezoutian. Since
$w(T^{-1})=-\eps w(T)$ the values of $w$ on $S^1$ are either real or
purely imaginary. Hence $w$ gives a continuous map $S^1\rightarrow
\bfP^1(\R)$ in either case. This yields a map $\Z\rightarrow \Z$ via
$H_1$ well defined up to sign. Choosing orientations appropriately
this map is again multiplication by the signature $\varsigma$ (defined
as zero in the skew symmetric case).

In practice one can compute $\varsigma$ using Sylverster's simple
characterization (that applies equally well to both the classical and
the skew cases). We associate to $w$ a word $\phi$ in two letters say
$A$ and $B$ as follows. Start with the empty word.  Traverse
$\bfP^1(R)$ or $S^1$ in the standard orientation starting at the base
point $\infty$ or $1$ respectively. Append $A$ (resp. $B$) to $\phi$
on the right if you encounter a root of $p$ (resp. $q$), including
multiplicities, finishing when you reach back the base point.  Now
recursively remove from $\phi$ any instance of repeated symbols $AA$
or $BB$. We end with a word consisting of $r$ pairs $\cdots ABAB
\cdots$ or $\cdots BABA \cdots$. Then~$\varsigma$ equals $r$ or $-r$
respectively.

In particular, $\Psi$ is definite if and only if the roots of $p$ and
$q$ interlace in the unit circle. This is one of the crucial
calculations of~\cite{BH} (see e. g. Theorem $4.8$ in \emph{loc. cit.}), which was done directly without
any reference to Hermite's result or its variants.

It is not hard to see~\cite[Prop. 2.3.3]{Ad} that the number of words corresponding
to signature $\varsigma$ is 
$$
\binom d{\tfrac12(d-\varsigma)}^2.
$$
(Necessarily $\varsigma\equiv d \bmod 2$; in fact in the symmetric
case, if $d$ is even then $\varsigma\equiv d \bmod 4$.) It follows
that we should expect the signature to be typically small if $p$ and
$q$ are picked in some random fashion. This appears to be indeed the
case. For example, considering all pair of coprime polynomials with
only cyclotomic factors and of degree $15$ with $\epsilon=-1$ we find
the following distribution of signatures
$$
\begin{array}{c|r|r|r|r|r|r|r|r|r}
\varsigma&15&13&11&9&7&5&3&1\\
\hline
\#&25& 118 &179 & 5935&  41242&  75458&  184173& 
268640
\end{array}
$$
with symmetrical values for $\varsigma=-1,-3,\ldots,-15$.

\subsection{Examples}
We end this section with some examples. The skew Bezoutian
construction can be done over a commutative ring (details will appear
in a later publication). Here we work over $\Z$.

Alternate constructions for the first two examples below can be found
in~\cite[Section $3$]{Bayer}. In loc. cit. the constructions use the
trace form (as mentioned in the introduction) with respect to an
extension $K/{\bf Q}$, where $K$ is a suitable cyclotomic
field. Examples in~\cite{Bayer} and~\cite[\S 4]{Bayer-Martinet} also include the Leech lattice, the
Coxeter--Todd lattice, etc. See also~\cite[\S 1]{Bayer-unimod} for related work 
where the question of the existence of a definite unimodular lattice with an isometry 
having a prescribed cyclotomic characteristic polynomial is addressed.

1) Let
$$
p=\Phi_1\Phi_2\Phi_3\Phi_5=x^8 + 2x^7 + 2x^6 + x^5 - x^3 - 2x^2 -
2x - 1,
\quad 
q=\Phi_{30}=x^8 + x^7 - x^5 - x^4 - x^3 + x + 1,
$$
where $\Phi_n$ is the $n$-th cyclotomic polynomial. Then
$$
w=-p/q=1 + x + x^2 + x^3 + x^4 + x^5 -x^{10}+ O(x^{11})
$$
and 
$$
B^*(p/q)=
\left(\begin{matrix}
 {2}&{1}&{1}&{1}&{1}&{1}&{0}&{0}\\
 {1}&{2}&{1}&{1}&{1}&{1}&{1}&{0}\\
 {1}&{1}&{2}&{1}&{1}&{1}&{1}&{1}\\
 {1}&{1}&{1}&{2}&{1}&{1}&{1}&{1}\\
 {1}&{1}&{1}&{1}&{2}&{1}&{1}&{1}\\
 {1}&{1}&{1}&{1}&{1}&{2}&{1}&{1}\\
 {0}&{1}&{1}&{1}&{1}&{1}&{2}&{1}\\
 {0}&{0}&{1}&{1}&{1}&{1}&{1}&{2}
\end{matrix}\right)
$$
The lattice $\Z[x]/(q)$ with this quadratic form is the well-known
$E_8$ lattice and $\gamma$ is a Coxeter element of the corresponding Weyl group.

2) Similarly the $A_n$ lattice with Cartan matrix 
$$
C_n:=\left(\begin{matrix}
2&-1&0&\cdots&&0&0\\
-1&2&-1&\cdots&&&0\\
& && \vdots&&&\\
0&&&\cdots&-1&2&-1\\
0&0&&\cdots&&-1&2
\end{matrix}\right)
$$
arises as the skew Bezoutian $B^*(p/q)$, where
$$
p= x^n-1 \qquad q=x^n+x^{n-1}+\cdots +x +1
$$
and $\gamma$ represents an $n$-cycle in $S_n$.

3) Let $q=x^{10} + x^9 - x^7 - x^6 - x^5 - x^4 - x^3 + x + 1$ be the
Lehmer polynomial (the integer polynomial of smallest known Mahler
measure bigger than $1$). For $p$, we search among the polynomials of
degree $10$ which are products of cyclotomics. We find eight such that
$B^*(p,q)$ is isometric to the unimodular lattice $I_{9,1}$ of signature
$(9,1)$. These are tabulated below.
$$
\begin{array}{ll}
\Phi_1^3\Phi_2\Phi_3\Phi_5 & x^{10} - x^8 - x^7 + x^3 + x^2 - 1\\
\Phi_1\Phi_2^3\Phi_3\Phi_5 & x^{10} + 4x^9 + 7x^8 + 7x^7 + 4x^6 - 4x^4
- 7x^3 - 7x^2 - 4x - 1 \\
\Phi_1\Phi_2\Phi_3\Phi_5\Phi_6 & x^{10} + x^9 + x^8 + x^7 + x^6 - x^4
- x^3 - x^2 - x - 1 \\
\Phi_1\Phi_2\Phi_3\Phi_7 & x^{10} + 2x^9 + 2x^8 + x^7 - x^3 - 2x^2 -
2x - 1 \\
\Phi_1\Phi_2\Phi_3\Phi_9 & x^{10} + x^9 - x - 1\\
\Phi_1\Phi_2\Phi_3\Phi_{18}& x^{10} + x^9 - 2x^7 - 2x^6 + 2x^4 + 2x^3
- x - 1 \\
\Phi_1\Phi_2\Phi_5\Phi_8 & x^{10} + x^9 + x^6 - x^4 - x - 1 \\
\Phi_1\Phi_2\Phi_5\Phi_{10} & x^{10} - 1
\end{array}
$$
We do not know if these isometries are in the same conjugacy class.

4) In the paper~\cite{MPV} the authors consider the modification
$q(x):=p(x)\pm x^m$ of a monic reciprocal polynomial $p$ of even
degree $2m$ consisting of adding a single monomial $\pm x^m$. The skew
Bezoutian $B^*(p,q)$ then yields a skew-symmetric form of determinant
$\Res(p,q)=1$ and a symplectic transformation of characteristic
polynomial $q$. For example, if we again take $q$ to be the Lehmer
polynomial we see that it is also the characteristic polynomial of a
symplectic transformation. As pointed out in~\cite[\S 4]{MPV} it is
remarkable that $q(x)+x^5$ is actually a product of cyclotomic
polynomials.

In light of Theorem~\ref{characterization}, the modification used
in~\cite{MPV} can be seen as an example of modifying a symplectic
transformation by multiplying it by a single transvection. It would be
interesting to extend their results and study how this modification
affects the Mahler measure of the characteristic polynomial.

\section{Isometries with given characteristic polynomial}

The goal of this section is to give a new and effective proof of the following 
well-known result (see, e.g.~\cite[Lemma $1.2$ and remarks $1.3$, $1.4$, $1.5$]{Mil}),
 using the skew Bezoutian. 

  We keep the notation of \S\ref{skew-bezoutian}. 

  \begin{theorem}
\label{prop-answer}
Let $q\in k[T]$ be a monic reciprocal polynomial of degree $d\geq
1$. Then

    1) There exists a non-degenerate symmetric bilinear space over $k$
    of dimension $d$ with an isometry of characteristic
    polynomial $q$.

    2) If, in addition, $d$ is even there exists a non-degenerate
    skew-symmetric bilinear space over $k$ of dimension $d$
     with an isometry of characteristic
    polynomial $q$.
  \end{theorem}
  
  For related work where a quadratic structure is prescribed as well see~\cite{Bayer1}.
  
\begin{proof}
  The main idea is to use the skew Bezoutian. If we can find a
  polynomial $p\in k[T]$, which is $(-\eps)$-reciprocal and coprime to
  $q$ then the skew Bezoutian $B^*(p,q)$ provides an explicit answer to
  what we are looking for. As discussed above, the skew Bezoutian
  comes equipped with an isometry of characteristic polynomial $q$ and
  is non-degenerate if $p$ and $q$ are coprime. Knowledge of
  $\Res(p,q)$ will help us show that the bilinear form we construct is non-degenerate by
  Proposition~\ref{bez-res}. 

  In the skew-symmetric case, where $\eps=-1$ and $d$ is
  assumed even, we may always find such a $p$. Indeed, the polynomial
$$
   p(T):=q(T)+T^{m},
$$
where $m:=d/2$ satisfies all the requirements we need: $p$ is clearly
reciprocal and coprime to $q$. Moreover we may easily compute $\Res(p,q)$. It is $(\prod_b b)^{m}$,
where $b$ runs over the roots of $q$ counted with multiplicity, and
this equals $q(0)^m$. Since $q$ is monic and reciprocal
$q(0)=1$. See the remark following the proof for a discussion on the relevance of this computation.


   Now we turn to the case where $\eps=1$. Let $\mathcal Q_0$ be the polynomial of $k[T]$ 
   such that one has the factorisation
\begin{equation}\label{defQ0}
     q(T)=(T-1)^{v_+}(T+1)^{v_-}\mathcal Q_0(T),\qquad \mathcal Q_0(\pm 1)\not=0\,.
\end{equation}
     
As $q$ is reciprocal, by the observations of \S\ref{prelim}, its
order of vanishing $v_+$ at $1$ is even. Hence, $\mathcal Q_0$ is also
reciprocal; let $d_0$ be its degree.  Assume for the moment that
$d_0>0$ and set
$$
     \mathcal P_0(T):=(T-1)^{e}(T+1)^{d_0-e}\,,
$$
for some odd  integer $0\leq e\leq d_0$. Note that $d_0-e$
is odd also since $d_0$ is even as $\mathcal Q_0$ is a reciprocal polynomial
not vanishing at $-1$ (see \S\ref{prelim}).
     
By construction $\mathcal P_0$ is monic, skew-reciprocal, of degree $d_0$ and
coprime to $\mathcal Q_0$. The skew Bezoutian $(V_0,\Psi_0)$ of $\mathcal P_0,\mathcal Q_0$ is
then a non-degenerate symmetric bilinear space over $k$ of dimension
$d_0$. The corresponding isometry $\gamma_0$ has characteristic
polynomial $\mathcal Q_0$. To obtain the space $V$ we are after we consider
$$
V:=V_0\perp V_+\perp V_-
$$
where $V_{\pm}$ is a vector space over $k$ of dimension $v_\pm$. We
put on $V_\pm$ an arbitrary non-degenerate symmetric bilinear form
$\Psi_\pm$ and consider $\Psi:=\Psi_0\perp \Psi_+\perp\Psi_-$
and $\gamma:=\gamma_0\perp \id_{V_+}\perp (-\id_{V_-})$. It is now
clear that $(V,\Psi)$ and $\gamma$ fulfill the requirements.

The same construction works if $d_0=0$; just ignore $V_0$
altogether. This completes the proof of 1).  
\end{proof}

\begin{remark}\label{DetOverRings}
Note that the proof actually gives a (skew-)symmetric space 
and an isometry of characteristic polynomial~$q$ defined over the ring 
of coefficients of the polynomials~$p$ and~$q$. 
In case $2)$ of Theorem~\ref{prop-answer} the determinant of this space is $1$ 
as can be seen from the computation performed in the proof. 
For case $1)$, see Remark~\ref{DetOverRings2}.
\end{remark}

{} It
seems natural to try and compute other invariants attached to the bilinear
space constructed in terms of the polynomials $p$ and $q$. In the
following section we focus on the case $\eps=1$ and
 we investigate how the spinor norm of the isometry
constructed in the proof of Theorem~\ref{prop-answer} can
be expressed in terms of the polynomial~$q$.
  
\section{Spinor norm of an isometry with prescribed characteristic
  polynomial}\label{spin-prescribed}
  
Recall that if $(V,\Psi)$ is a non-degenerate finite dimensional quadratic space, the
spinor norm of an isometry of $(V,\Psi)$ can be defined as follows:
first let $v$ be a non isotropic vector of $V$ and let $r_v$ be the
reflection with respect to the hyperplane $v^\bot$. We define the
spinor norm $\nsp(r_v)$ to be the class in $k^\star/(k^\star)^2$ of
$\Psi(v,v)$. Now any isometry $\sigma$ of $V$ is a product $\prod_v
r_v$, where $v$ runs over a finite set of non-isotropic vectors of
$V$. It is known that $\sigma\mapsto \prod_v \Psi(v,v)$ gives a
well-defined {\it spinor norm} homorphism
  $$
 \nsp:  O(V,\Psi)\lra k^\star/(k^\star)^2\,,
  $$
  which is onto as soon as $d:=\dim V\geqslant 2$. Note in particular
  that $\nsp(-\id_V)=\det(V,\Psi)$, where $\det(V,\Psi):=\det\left(
    \Psi(v_i,v_j)\right)$ for any basis $v_1,\ldots,v_d$ of $V$. (If
  $v_1,\ldots,v_d$ is an orthogonal basis of $V$ then
  $-\id_V=\prod_{i=1}^d r_{v_i}$ and $\det V=\prod_{i=1}^d
  \Psi(v_i,v_i)$ .)

  We recall the following formula due to Zassenhaus
  (see~\cite[p. 444]{Za}) which gives a useful way to compute the
  spinor norm of an isometry. To state and prove the results of this 
  section it will be convenient to use the following notation introduced by 
  Zassenhaus in his original paper. If $\sigma$ is an endomorphism 
  of $V$ and if $\lambda\in k$ then we let $M(\lambda,\sigma)$ 
  be the maximal subspace of $V$ on which $\sigma-\lambda{\id_V}$ acts as a 
  nilpotent endomorphism of $V$. In particular the dimension  of $M(\lambda,\sigma)$ 
  is the multiplicity of $\lambda$ as a root of the characteristic polynomial of $\sigma$.
  
  
  \begin{theorem}[Zassenhaus]\label{Z-formula}
  Let $\gamma$ be an isometry of a non-degenerate quadratic space
  $(V,\Psi)$ over $k$ and let $v_{\pm}$ be the dimension of $M(\pm 1,\gamma)$. 
Let $q$ be the characteristic polynomial of $\gamma$. Then $M(-1,\gamma)$ is non-degenerate and,
  if we denote by $q_-$ the polynomial such  that
$$
q(T)=(T+1)^{v_-}q_-(T)\,,\qquad q_-(-1)\not=0\,,
$$
then
$$
\nsp(\gamma)=\det\left(M(-1,\gamma),\Psi\right)(-2)^{-(\dim V-v_-)}q_-(-1)\,,
$$
in $k^\star/(k^\star)^2$.
\end{theorem}

  
 \begin{proof}
  Let us describe the main ideas of the proof based on 
  Zassenhaus original paper~\cite[pp. 444--446]{Za}.
  Let us consider the subspace of $V$:
  $$
  \widehat{M}(-1,\gamma):=\bigcap_{n\geq 1} \left(\sigma+{\id_V}\right)^n V\,.
  $$
  Then Zassenhaus shows~(\cite[Prop. 2 p. 437 \& its corollary p. 438]{Za}) that one has the orthogonal splitting
  $$
  V=M(-1,\gamma)\perp\widehat{M}(-1,\gamma)\,,
  $$
  thus both these spaces are non-degenerate with respect to the restriction of $\Psi$.
In particular the formula
\begin{equation}\label{Z-nsp}
{\rm sn}(\gamma):=\det\left(M(-1,\gamma),\Psi\right)\cdot \det\left(\frac{\gamma+{\id_V}}{2}\mid \widehat{M}(-1,\gamma)\right)\,,
\end{equation}
defines a function on the orthogonal group $O(V,\Psi)$ with values in the classes modulo non zero squares of $k^\star$. 
Zassenhaus then shows~(\cite[Theorem p. 446]{Za}) that the map {\rm sn} is a group homomorphism and that it 
coincides with $\nsp$ (see~\cite[(2.10b) p. 446]{Za}).


 One has $M(-1,\gamma)^\perp=\widehat{M}(-1,\gamma)$ and the restriction of $\gamma$ to $M(-1,\gamma)^\perp$ has characteristic polynomial $q_-$. Therefore~\eqref{Z-nsp} yields
$$
\nsp(\gamma)=\det\left(M(-1,\gamma),\Psi\right)\cdot (-2)^{-(\dim V-v_-)}q_-(-1)\,,
$$
in $k^\star/(k^\star)^2$, which completes the proof.
\end{proof}

  \begin{corollary} \label{nspin-formula}
 With notation as above fix an isometry $\gamma$ of $(V,\Psi)$. Let $q\in k[T]$ be the characteristic 
 polynomial of $\gamma$ and let $\mathcal Q_0\in k[T]$ be as in~\eqref{defQ0}.
 Then the spinor norm of $\gamma$
  is given by  
  $$
  \nsp(\gamma)= \mathcal Q_0(-1)\det\left(M(-1,\gamma),\Psi\right)\,,
  $$ 
  in $k^\star/(k^\star)^2$.
  \end{corollary}

  \begin{proof}
With notation of Theorem~\ref{Z-formula} one has
$$
q_-(T)=(T-1)^{v_+}\mathcal Q_0(T)\,.
$$

 We deduce
$$
\nsp(\gamma)=\det\left(M(-1,\gamma),\Psi\right)(-2)^{-(\dim V-v_-)}(-2)^{v_+}\mathcal Q_0(-1)\,,
$$
in $k^\star/(k^\star)^2$.

Therefore:
\begin{align*}
\nsp(\gamma)&=\det\left(M(-1,\gamma),\Psi\right)(-2)^{-(\dim V-(v_-+v_+))}\mathcal Q_0(-1)\\
&=\det\left(M(-1,\gamma),\Psi\right)(-2)^{d_0}\mathcal Q_0(-1)\,,
\end{align*}
modulo nonzero squares. That is the desired formula since $d_0:=\deg \mathcal Q_0$ is even.
  \end{proof}

 From the above corollary we further deduce how to decide when we can  prescribe
 the spinor norm  and the characteristic polynomial of an isometry. 
 \par

\begin{corollary}
\label{cor-spin}
Let $q\in k[T]$ be a monic reciprocal polynomial of degree $d\geq
1$ and let $\mathcal Q_0\in k[T]$ be as in~\eqref{defQ0}.

(i) If $v_-(q)>0$ then there exists a non-degenerate symmetric
bilinear space over $k$ of dimension $d$ with an isometry $\gamma$ of
characteristic polynomial $q$ and arbitrary spinor norm
$\nsp(\gamma)$.  In particular, this is true if $d$ is odd.
\par
(ii) If $v_-(q)=0$  and  $\gamma$ is an isometry with characteristic
polynomial $q$ then its spinor norm equals $\mathcal Q_0(-1)$ (modulo nonzero
squares). In particular, this is the case if $q$ is separable and $d$
is even.  
\end{corollary}
\begin{proof}
  (i) Fix a representative $s$ for a class in $k^\star/(k^\star)^2$. If $v_->0$ 
  we can always choose $V_-$ to have $\det(V_-)\equiv
  s\mathcal Q_0(-1) \bmod (k^\star)^2$.  The result now follows from
  Corollary~\ref{nspin-formula}. If $d$ is odd by the observations
  of \S\ref{prelim} $v_-$ is odd and hence positive.
\par

(ii) The first statement follows from Corollary~\ref{nspin-formula}.
 Assume $q$ to be separable; if $v_-(q)>0$ then the quotient
 $q(T)/(T+1)$ is a reciprocal polynomial of odd degree. So $-1$ is
 also a root of the quotient which contradicts the separability of
 $q$. 
\end{proof}
  
\section{Discriminant of a quadratic space having an isometry with prescribed characteristic polynomial}   
 This section is devoted to the study of the relation between the discriminant of a quadratic space $(V,\Psi)$
 and the characteristic polynomial of an isometry of $O(V,\Psi)$. If $(V,\Psi)$ is a
quadratic space over $k$ we let its {\it discriminant} be $\disc
(V,\Psi):=(-1)^{n(n-1)/2}\det (V,\Psi)$ where $d:=\dim V$. 

The results we present here are well-known. 
 The idea emphasized in the following statement (that can be found, e.g., in~\cite[Th. 3.4]{Mil}) is that to an $\varepsilon$-symmetric 
non-degenerate bilinear space $(V,\Psi)$ equipped with an isometry $\gamma$, 
we can naturally associate a $(-\varepsilon)$-symmetric
 non-degenerate bilinear space $(V,\Psi_\gamma)$. 

 \begin{lemma}\label{lemma-star}
 Let $(V,\Psi)$ be an $\eps$-symmetric non-degenerate bilinear space and let $\gamma$ be an isometry of 
$(V,\Psi)$. 
We define the bilinear form $\Psi_\gamma$ on $V$ by:
$$
\Psi_\gamma(u,v)=\Psi\left((\gamma-\gamma^{-1})(u),v\right),\qquad u,v\in V\,.
$$

Denoting as before by $q$ the characteristic polynomial of $\gamma$, we have:

 $(i)$ $(V,\Psi_\gamma)$ is $(-\varepsilon)$-symmetric, 

$(ii)$ $\det(V,\Psi_\gamma)=q(1)q(-1)\det\gamma \det(V,\Psi)$, 

$(iii)$ $\gamma$ is an isometry of the bilinear space $(V,\Psi_\gamma)$. 
\end{lemma}

\begin{proof}
First note that for any isometry $\gamma$ of a bilinear space $(V,\langle\cdot,\cdot\rangle)$ and any element $h\in k[x,x^{-1}]$ we have
\begin{equation} \label{isom-bil}
\langle h(\gamma)u,v\rangle = \langle u,h(\gamma^{-1})v\rangle, \qquad u,v \in
V.
\end{equation}

For $(i)$ we fix $u,v\in V$ and we compute, using~(\ref{isom-bil}),
$$ 
\Psi_\gamma(v,u)=\Psi\left((\gamma-\gamma^{-1})v,u\right)
                =\Psi\left(v,(\gamma^{-1}-\gamma)u\right)
                =-\Psi\left(v,(\gamma-\gamma^{-1})u\right)\,.
$$

The right hand side equals $-\eps\Psi_\gamma(u,v)$ since $\Psi$ is $\eps$-symmetric.

For $(ii)$, we denote by $d$ the dimension of $V$ and we fix a basis
$\mathcal{B}=(e_1,\ldots,e_d)$ of $V$. Let $\mathcal{Q}$ be the Gram
matrix of $\Psi$ with respect to $\mathcal{B}$ and $M$ be the matrix
representation of $\gamma$ in the basis $\mathcal{B}$.  The Gram
matrix of $\Psi_\gamma$ with respect to $\mathcal{B}$ is
$\left(\Psi\left((\gamma-\gamma^{-1})e_i,e_j\right)\right)_{i,j}$.
That matrix equals ${}^t\left((M-M^{-1})\right)\mathcal{Q}$. Taking
determinants we get
$$
\det(V,\Psi_\gamma)=\det(M-M^{-1})\det(V,\Psi)=\det(\gamma-\gamma^{-1})\det(V,\Psi)\,,
$$
which is the formula we wanted since $\det(\gamma-\gamma^{-1})=\det(\gamma)\det(\gamma^2-\id_V)=\det(\gamma)q(-1)q(1)$.

Finally, $(iii)$ is a straightforward consequence of the fact that
$\gamma$ and $\gamma-\gamma^{-1}$ commute.
\end{proof}

The construction of the bilinear form $\Psi_\gamma$ from the data $(\Psi,\gamma)$
 can be iterated to produce a sequence of bilinear forms $\Psi_0=\Psi$,
  $\Psi_1=\Psi_\gamma$, and more generally for any $j\geq 0$:
  $$
  \Psi_j\colon (u,v)\in V\times V\mapsto \Psi\left((\gamma-\gamma^{-1})^{j}u,v\right)\,.
    $$
The generalization of Lemma~\ref{lemma-star} to $\Psi_j$ is straightforward and can be found 
in~\cite[Th. 3.4]{Mil}. Using Lemma~\ref{lemma-star} we deduce the following statement.

\begin{proposition}
\label{det-isom}
Let $q\in k[x]$ be a monic reciprocal polynomial with $q(\pm 1)\neq
0$. Then the discriminant $\disc(V,\Psi)$ of a non-degenerate
quadratic space $(V,\Psi)$ over $k$ with an isometry of characteristic
polynomial $q$ is uniquely determined.  More precisely, for any such
space we have
$$
  \det(V,\Psi) \equiv q(-1)q(1) \bmod (k^\star)^2.
$$
\end{proposition}

This statement is well-known and can be found e.g. in~\cite[\S 7, Lemma c)]{Le}.
 
 \begin{proof}
  We invoke Lemma~\ref{lemma-star}$(ii)$ in the case $\eps=1$. 
 Indeed $\det \gamma=1$ since $q$ is reciprocal. Moreover the formula for $\det(V,\Psi_\gamma)$ 
 implies that $\Psi_\gamma$ is non-degenerate by our assumption on $q$ and $\Psi$.
Thus $\det(V,\Psi_\gamma)$
 is a square since $(V,\Psi_\gamma)$ is a non-degenerate 
 skew-symmetric bilinear space.
 \end{proof}
 
 \begin{remark}\label{DetOverRings2}
 In fact, from the proof of Theorem~\ref{prop-answer} we see that for every 
 odd integer $e$ in the range $0\leq e\leq \deg q=d$ we can find $(V,\Psi)$ defined 
 over the ring of coefficients of $q$ that satisfies:
 $$
 \det(V,\Psi)=q(1)q(-1)^{d-e}\,.
 $$
 \end{remark}
 
From the above proposition we deduce the following corollary that answers the question  
investigated in this section. 

\begin{corollary}\label{cor-BE}
Let $(V,\Psi)$ be a
     non-degenerate quadratic space over $k$.     
      
     $(i)$ Let $q$ be a
     reciprocal polynomial which is the characteristic polynomial of
     an isometry $\gamma$ of $(V,\Psi)$. 
Then with notation as in~\S\ref{spin-prescribed} (in particular we use the factorisation~\eqref{defQ0}),
       \begin{equation}
\label{disc-fmla}
   \det(V,\Psi)\equiv \det\left (M(-1,\gamma)\right)\det\left(M(1,\gamma)\right)\mathcal Q_0(-1)\mathcal Q_0(1)\bmod
   (k^\star)^2. 
 \end{equation}
 
 $(ii)$ Let $q$ be a separable
  reciprocal polynomial in $k[T]$ of even degree.  If there exists an
  isometry $\gamma$ of a non-degenerate quadratic $k$-space $(V,\Psi)$
  of characteristic polynomial $q$ then $\disc(V,\Psi)\equiv\disc(q)
  \bmod (k^\star)^2$.

\end{corollary}

\begin{proof}
 From~\cite[Prop. 2 \& its corollary]{Za} and since we assume ${\rm char}\, k\neq 2$, one easily deduces the orthogonal decomposition:
 \begin{equation}\label{orth-split}
 V=M(1,\gamma)\perp M(-1,\gamma)\perp\left(\widehat{M}(1,\gamma)\cap\widehat{M}(-1,\gamma)\right)\,.
 \end{equation}
 Thus
 $$
 \det V=\det M(1,\gamma)\det M(-1,\gamma)\det \left(\widehat{M}(1,\gamma)\cap\widehat{M}(-1,\gamma)\right)\,,
 $$
 where the quadratic structure on each vector space is given by the suitable restriction of $\Psi$.
 
  Each subspace on the right hand side of~\eqref{orth-split} is stable under $\gamma$ and by defintion 
  of the subspaces $M(\pm 1,\gamma)$, the restriction of $\gamma$ to $\widehat{M}(1,\gamma)\cap\widehat{M}(-1,\gamma)$ 
  has characteristic polynomial $\mathcal Q_0$. Thus $(i)$ follows by applying Proposition~\ref{det-isom}.

The statement $(ii)$ is an easy consequence of $(i)$ and the following
well-known lemma.
\end{proof}

\begin{lemma} \label{lemme-edwards}
  Let $q\in k[x]$ be a monic separable reciprocal polynomial of even degree
  $2m$. Then
$$
\disc q \equiv (-1)^mq(-1)q(1) \bmod (k^\star)^2.
$$
\end{lemma}

\begin{proof}
  The hypothesis on $q$ guarantees that $q(\pm 1)\neq 0$, i.e.,
  $q=\mathcal Q_0$. Indeed, if $q$ is reciprocal then $v_+$ must be even. If in
  addition $q$ is separable then $v_+(q)=0$. As we argued in the proof
  of Corollary~\ref{cor-spin} (ii) we also have $v_-(q)=0$.

  We may assume without loss of generality that $q$ is
  irreducible. Let $K:=k[x]/(q)$. The extension $K/k$ is separable and
  $\disc q$ is the discriminant of the quadratic space $(K,\Psi)$,
  where $\Psi(a,b):=\Tr_{K/k}(ab)$.  A calculation like that in
  \cite[Prop. A.3]{McG} (see also the discussion at the beginning of section $2$ in~\cite{Ba}) 
  finishes the proof. (Let $L\subseteq K$ be the
  subfield fixed by the involution $x\mapsto x^{-1}$, and let $\N_{L/k}$ denote the norm map relative to 
  $L/k$. Then
  $K=L(x-x^{-1})$. The subspaces $L$ and
  $(x-x^{-1}) L$ are orthogonal hence $\det K =
  \N_{L/k}(x-x^{-1})\,\det L^2$ and $\N_{L/k}(x-x^{-1})=q(-1)q(1)$.)
\end{proof}

For an alternate proof of the lemma  see~\cite[proof of Th. 2]{E}. The statement $(ii)$ of Corollary~\ref{cor-BE} 
can be found e. g. in~\cite[Th. (1.2)]{Ba}.
 
 \section{Isometries with given Jordan form}\label{sectionjordan}

 We end with a characterization of the Jordan form of isometries of
 non-degenerate bilinear spaces. The main result goes back to (at
 least) Wall \cite{Wall} (see also~\cite{HM}, \cite[section 3]{Mil} and~\cite[IV, 2.15 (iii)]{SpSt}). We include a proof for
 the reader's convenience using the skew Bezoutian to construct the
 isometries.

 We assume our field $k$ is now algebraically closed (and of
 characteristic different from~$2$ as before). Fix a vector space $V$
 of dimension $r$ over $k$. For $\gamma\in \End(V),\lambda\in k^\star$
 and $m\in \N$, let $\mu(\gamma;\lambda,m)$ be the number of Jordan blocks
 of $\gamma$ of size $m$ and eigenvalue $\lambda$.

 We start with a few preparatory results. 
 The following crucial statement, (very close to the first part of~\cite[Th. 3.2]{Mil}), 
 will help us perform a reduction step needed in the proof of Theorem~\ref{jordanblocks}.

\begin{lemma}\label{th32-mil}
 Let $(V,\Psi)$ be a non-degenerate $\eps$-symmetric space
  equipped with a unipotent isometry $\gamma$. 
  We have an orthogonal splitting:
$$
V=\perp_{m\geqslant 1} V ^{(m)}\,,
$$  
where $\gamma$ acts on each $V^{(m)}$ as a sum of Jordan blocks $J_1(m)$.
 In particular, each $(V^{(m)},\Psi)$ is a non-degenerate $\eps$-symmetric space.
\end{lemma}

\begin{proof}
 Let $n$ be the largest index $m$ with $V^{(m)}\neq 0$. We
claim that $V^{(n)}$ is non-degenerate.

Since $\gamma$ is unipotent and preserves $\rad
(V^{(n)},\Psi)$ we have $\rad (V^{(n)},\Psi)\subseteq \ker (\gamma-{\id_V})=\im
\left((\gamma-{\id_V})^{n-1}\right)$.


 Taking $h(x)=(x-1)^{n-1}$ in
\eqref{isom-bil} it follows that $\rad (V^{(n)},\Psi)\subseteq \rad (V,\Psi)$
proving our claim. We deduce that $V^{(n)}$ splits off from $V$ as an orthogonal direct summand.

 We conclude by finite descending induction on $m\geqslant 1$.
\end{proof}

The following Lemma can be seen as a complement to Lemma~\ref{lemma-star}. 
In the notation of Lemma~\ref{lemma-star} it gives an additional property of $\Psi_\gamma$ in
the case where $-1$ is not an eigenvalue of $\gamma$. For any bilinear
space $(W,\langle\cdot,\cdot\rangle)$, its \emph{radical} ${\rm
  rad}(W)$ is the subspace $\{v\in W\colon \langle v,w\rangle=0\text{
  for all }w\in W\}$.

\begin{lemma}\label{lemma-N}
 With notation as in Lemma~\ref{lemma-star}, we assume further that $\gamma+{\id_V}$ is invertible. 
 Then we have
$$
{\rm rad}(V,\Psi_\gamma)=\ker(\gamma-{\id_V})\,.
$$
\end{lemma}

\begin{proof}
 Fix a vector $u\in V$. We have $u\in {\rm rad}(V,\Psi_\gamma)$ if and only if
  $\Psi\left((\gamma-\gamma^{-1}) u,v\right)=0$ for all $v\in V$. Since $\Psi$
   is non-degenerate, this is equivalent to $(\gamma-\gamma^{-1}) u=0$, i.e.
    $(\gamma^2-{\id_V})u=0$. Rewriting the last equation 
$$
(\gamma+{\id_V})\circ (\gamma-{\id_V})u=0\,,
$$
the lemma follows since we have assumed $\gamma+{\id_V}$ to be invertible. 
\end{proof}


\begin{remark}\label{iteration}
 As for the case of Lemma~\ref{lemma-star} the generalization of Lemma~\ref{lemma-N} to $\Psi_j$ 
 is straightforward. Let us mention for example that if $\Psi$ is $\eps$-symmetric then $\Psi_j$ is
  $(-1)^j\eps$-symmetric with radical $\ker\left((\gamma-{\id_V})^j\right)$ (see~\cite[Th. 3.4]{Mil} where 
  the general version of the construction is used).
\end{remark}

\begin{corollary}\label{cor-lemma-N}
 With hypotheses as in Lemma~\ref{lemma-star}, assume $\eps=1$ and $\gamma$ unipotent.
  Consider the Jordan block decomposition of $\gamma$:
$$
\bigoplus_{i=1}^r J_{m_i}(1)\,, \qquad m_1\leqslant m_2\leqslant\cdots\leqslant m_r,\, \sum_i{m_i}=\dim V\,,
$$
where $J_{m_i}(1)$ stands for the Jordan block of size $m_i$ attached to the eigenvalue $1$. We have

$$
\sum_{i=1}^r{(m_i-1)}\equiv 0\,(\mathrm{mod}\, 2)\,.
$$

 In particular there are evenly many indices $i$ for which $m_i$ is even. 

\end{corollary}

\begin{proof}
 Since $\left(V/{\rm rad}(V,\Psi_\gamma),\Psi_\gamma\right)$ is a non-degenerate
  skew-symmetric space, its dimension is even. It follows then from Lemma~\ref{lemma-N} that $\dim V\equiv 
\dim \ker(\gamma-{\id_V})\,(\mathrm{mod}\,2)$.
Since $\gamma$ is unipotent its number of Jordan blocks equals $\dim \ker(\gamma-{\id_V})$ therefore
$$
\sum_{i=1}^r{m_i}\equiv r\,(\mathrm{mod}\,2)\,.
$$ 

 Equivalently $\sum_{i=1}^r{(m_i-1)}$ is even.
 
\end{proof}

We can now state and prove the main result of this section.

\begin{theorem} \label{jordanblocks}
Let $\gamma\in \End(V)$. Then $\gamma$ preserves a non-degenerate
    $\eps$-symmetric bilinear form on $V$ if and only if 

(i)
$$
\mu(\gamma,\lambda,m)=\mu(\gamma,\lambda^{-1},m), \qquad \lambda\neq
\pm 1, \quad  m\in \N,
$$
and

(ii)
$$
 (m-\delta)\mu(\gamma,\pm 1,m) \equiv 0 \bmod 2, \qquad m \in \N,
$$
where $\delta:= \tfrac12(1+\eps)$.
\end{theorem}
  
 \begin{proof}
   We give details for the orthogonal case $\eps=1$ the symplectic
   case $\eps=-1$ is completely analogous. For $m\geqslant 1$ let
   $J_m(\lambda)$ denote the Jordan block with size $m$ and eigenvalue
   $\lambda$.

   First we exhibit an isometry with a prescribed Jordan form
   satisfying the hypothesis (i) and (ii). Identify $V$ with $k^d$. If 
   $\gamma\in \End(V)$ is an endomorphism having Jordan form
   $M=J_m(\lambda)\oplus J_m(\lambda^{-1})$ with $\lambda\neq
   \lambda^{-1}$ consider $q=(T-\lambda)^m(T-\lambda^{-1})^m$. By
   Theorem~\ref{prop-answer} there exists a skew-reciprocal polynomial
   $p\in k[T]$ such that $q$ is the characteristic polynomial of an
   isometry of the non-degenerate quadratic space determined by
   $B^*(p,q)$, which by Theorem~\ref{characterization} has Jordan form
   $M$.  A similar argument applies to $J_m(\pm 1)$ for $m$ odd taking
   $q=(T-1)^m$ and $p=(T+ 1)^m$.

   Finally, let $m$ be even and set again $p:=(T+1)^m$ and
   $q:=(T-1)^m$. Now $U:=B^*(p,q)$, however, is skew-symmetric.
   Consider instead the symmetric matrix
$$
 A=  \left(
\begin{array}{cc}
0 & U \\
-U & 0
\end{array}
\right).
$$
 Since $p$ and $q$ are relatively prime $U$ and
hence also $A$ yield non-degenerate bilinear pairings. By
Theorem~\ref{characterization} and Theorem~\ref{prop-answer} there
exists $\gamma^\pm$ with Jordan form $J_m(\pm1)$ preserving
$U$. The map $\gamma:=\gamma^\pm\oplus \gamma^\pm$ then preserves $A$
giving our desired isometry.

We now show that the conditions on the multiplicities of the Jordan
blocks are necessary. Suppose then that $\gamma\in \End(V)$ preserves a
non-degenerate, symmetric bilinear pairing $\langle
\cdot,\cdot\rangle$ on $V$.
It follows that as $k[x,x^{-1}]$-modules $V^*\simeq V$. This implies
(i). 

For $\lambda\in k^\star$ let $V_\lambda\subseteq V$ be the subspace
annihilated by some power of $\gamma-\lambda$ and let
$W_\lambda:=V_\lambda\oplus V_{\lambda^{-1}}$ if $\lambda\neq
\lambda^{-1}$ and $W_{\pm 1}:=V_{\pm 1}$. Taking
$h(x)=(x-\lambda)(x-\lambda^{-1})$ or $x-(\pm1)$ in \eqref{isom-bil}
we see that the distinct non-zero $W_\lambda$'s are mutually
orthogonal with orthogonal sum $V$ and, in particular, they are
non-degenerate. To prove (ii) we may hence assume without loss of
generality that $\gamma$ is unipotent so $V=V_1$.

 Applying Lemma~\ref{th32-mil}, we can restrict further to the case where $V=V^{(m)}$ is a non-degenerate 
quadratic space on which $\gamma$ acts as a sum of $\mu(\gamma,1,m)$ Jordan blocks $J_1(m)$.




 Applying Corollary~\ref{cor-lemma-N} to $\gamma$, we deduce that $\mu(\gamma,1,m)(m-1)$ is even, which
is what we wanted to prove.

 
 Note that in the skew-symmetric case $(ii)$ follows directly from Lemma~\ref{th32-mil}. Indeed 
 $\gamma$ restricts to a unipotent isometry of the non-degenerate skew-symmetric space $(V^{(m)},\Psi)$. 
 Thus the dimension $m\mu(\gamma, 1,m)$ of this space is even.


\end{proof}

\end{document}